\documentclass{article}
\usepackage{abstract}
\usepackage{parskip}
\usepackage{amssymb}
\usepackage{amsmath}
\usepackage{amsthm}
\usepackage{tikz}
\usepackage{enumerate}
\usetikzlibrary{chains}
\newtheorem{theorem}{Theorem}[section]
\newtheorem{proposition}[theorem]{Proposition}
\newtheorem{lemma}[theorem]{Lemma}

\theoremstyle{definition}
\newtheorem{definition}[theorem]{Definition}
\theoremstyle{remark}
\newtheorem{remark}[theorem]{Remark}
\newtheorem{example}[theorem]{Example}
\newcommand{\chain}[4]{
        \filldraw (#2,1 - #3*#2) circle (0.05cm) node[below=5pt] {$C_{#4}(w)$};
        \foreach \y in {2,...,#1}
        {
            \draw (#2,\y - #3*#2 - 1.0) -- (#2,\y - #3*#2);
            \filldraw (#2,\y - #3*#2) circle (0.05cm);
        }
        }
\newcommand{\connectem}[5]{
        \foreach \y in {1,...,#1}
        {
            \draw (#2,#3 + \y) -- (#4,#5 + \y);
        }
        }
\newcommand{\bm}[1]{{\mbox{\boldmath $#1$}}}
\newcommand{\tand}{\text{ and }}
\newcommand{\bleq}{\leq_L}

\newcommand{\pleq}{\leq_S}

\newcommand{\minix}[3]{m_{#1,#2}(#3)}
\newcommand{\leh}[1]{\bm{c}(#1)}
\newcommand{\lehc}[2]{c_{#1}(#2)}
\newcommand{\lehm}[3]{c_{#1,#2}(#3)}

\newcommand{\chn}[2]{C_#1(#2)}
\newcommand{\inv}[1]{\text{\textnormal{Inv}}(#1)}
\newcommand{\ninv}[1]{\overline{\text{\textnormal{Inv}}}(#1)}

\newcommand{\card}[1]{|#1|}
\begin{document}
\title{A refinement of weak order intervals into distributive lattices}
\author{Hugh Denoncourt}
\maketitle
\begin{abstract}
In this paper we consider arbitrary intervals in the left weak
order on the symmetric group $S_n$. We show that
the Lehmer codes of permutations in an interval form a distributive
lattice under the product order. Furthermore, the rank-generating function of
this distributive lattice matches that of the weak order interval. We construct a poset such that its lattice of order ideals is isomorphic to the lattice of Lehmer codes 
of permutations in the given interval. We show that there are at least $\left(\lfloor\frac{n}{2}\rfloor\right)!$ permutations in $S_n$ that form a rank-symmetric 
interval in the weak order.
\end{abstract}
\section{Introduction and preliminaries}
\subsection{Introduction}
Our results concern intervals in the weak order of the symmetric group $S_n$. Intervals in this fundamental order can arise in unexpected contexts. For example, Bj{\"{o}}rner and Wachs \cite[Theorem 6.8]{permstatslinext} showed that the set of linear extensions of a regularly labeled two-dimensional poset forms an interval in the weak order. The Bell classes defined by Rey in \cite{bellorder} are also weak order intervals \cite[Theorem 4.1]{bellorder}.

Stembridge \cite[Theorem 2.2]{onfc} showed that the interval $\Lambda_w = [\text{id},w]$ in the weak order is a distributive lattice if and only if $w$ is a fully commutative element. The Lehmer code \cite{lehmer} is an $n$-tuple that encodes  information about the inversions of a permutation. Our main theorem, Theorem~\ref{t:distributivelattice}, states that the set of Lehmer codes for permutations in $\Lambda_w$, ordered by the product order on $\mathbb{N}^n$, is a distributive lattice. Furthermore, the rank-generating function of $\Lambda_w$ matches that of the corresponding distributive lattice. Theorem~\ref{t:distributivelattice} holds for arbitrary $w \in S_n$, so it tells us how an arbitrary weak order interval can be refined to form a distributive lattice when $w$ is not fully commutative.
\begin{center}
\begin{tikzpicture}[scale=1.0] 
\draw (-0.6,1.2) -- (0.0,0.0);
\draw (1.9,1.2) -- (0.0,0.0);
\draw (1.3,2.4) -- (-0.6,1.2);
\draw (1.3,2.4) -- (1.9,1.2);
\draw (0.7,3.6) -- (1.3,2.4);
\draw (2.5,2.4) -- (1.9,1.2);
\draw (1.9,3.6) -- (2.5,2.4);
\draw (1.3,4.8) -- (0.7,3.6);
\draw (1.3,4.8) -- (1.9,3.6);
\draw (3.8,2.4) -- (1.9,1.2);
\draw (3.2,3.6) -- (1.3,2.4);
\draw (3.2,3.6) -- (3.8,2.4);
\draw (2.6,4.8) -- (0.7,3.6);
\draw (2.6,4.8) -- (3.2,3.6);
\draw (4.4,3.6) -- (2.5,2.4);
\draw (4.4,3.6) -- (3.8,2.4);
\draw (3.8,4.8) -- (1.9,3.6);
\draw (3.8,4.8) -- (4.4,3.6);
\draw (3.2,6.0) -- (1.3,4.8);
\draw (3.2,6.0) -- (2.6,4.8);
\draw (3.2,6.0) -- (3.8,4.8);
\draw (0.0,0.0) node[fill=white] {$12345$};
\draw (-0.6,1.2) node[fill=white] {$21345$};
\draw (1.9,1.2) node[fill=white] {$12435$};
\draw (1.3,2.4) node[fill=white] {$21435$};
\draw (0.7,3.6) node[fill=white] {$31425$};
\draw (2.5,2.4) node[fill=white] {$13425$};
\draw (1.9,3.6) node[fill=white] {$23415$};
\draw (1.3,4.8) node[fill=white] {$32415$};
\draw (3.8,2.4) node[fill=white] {$12534$};
\draw (3.2,3.6) node[fill=white] {$21534$};
\draw (2.6,4.8) node[fill=white] {$31524$};
\draw (4.4,3.6) node[fill=white] {$13524$};
\draw (3.8,4.8) node[fill=white] {$23514$};
\draw (3.2,6.0) node[fill=white] {$32514$};
\draw (4.9,1.2) -- (5.5,0.0);
\draw (7.4,1.2) -- (5.5,0.0);
\draw (6.8,2.4) -- (4.9,1.2);
\draw (6.8,2.4) -- (7.4,1.2);
\draw (6.2,3.6) -- (6.8,2.4);
\draw (8.0,2.4) -- (7.4,1.2);
\draw (7.4,3.6) -- (6.8,2.4);
\draw (7.4,3.6) -- (8.0,2.4);
\draw (6.8,4.8) -- (6.2,3.6);
\draw (6.8,4.8) -- (7.4,3.6);
\draw (9.3,2.4) -- (7.4,1.2);
\draw (8.7,3.6) -- (6.8,2.4);
\draw (8.7,3.6) -- (9.3,2.4);
\draw (8.1,4.8) -- (6.2,3.6);
\draw (8.1,4.8) -- (8.7,3.6);
\draw (9.9,3.6) -- (8.0,2.4);
\draw (9.9,3.6) -- (9.3,2.4);
\draw (9.3,4.8) -- (7.4,3.6);
\draw (9.3,4.8) -- (8.7,3.6);
\draw (9.3,4.8) -- (9.9,3.6);
\draw (8.7,6.0) -- (6.8,4.8);
\draw (8.7,6.0) -- (8.1,4.8);
\draw (8.7,6.0) -- (9.3,4.8);
\draw (5.5,0.0) node[fill=white] {$00000$};
\draw (4.9,1.2) node[fill=white] {$10000$};
\draw (7.4,1.2) node[fill=white] {$00100$};
\draw (6.8,2.4) node[fill=white] {$10100$};
\draw (6.2,3.6) node[fill=white] {$20100$};
\draw (8.0,2.4) node[fill=white] {$01100$};
\draw (7.4,3.6) node[fill=white] {$11100$};
\draw (6.8,4.8) node[fill=white] {$21100$};
\draw (9.3,2.4) node[fill=white] {$00200$};
\draw (8.7,3.6) node[fill=white] {$10200$};
\draw (8.1,4.8) node[fill=white] {$20200$};
\draw (9.9,3.6) node[fill=white] {$01200$};
\draw (9.3,4.8) node[fill=white] {$11200$};
\draw (8.7,6.0) node[fill=white] {$21200$};
\end{tikzpicture}
\end{center}
\begin{center}
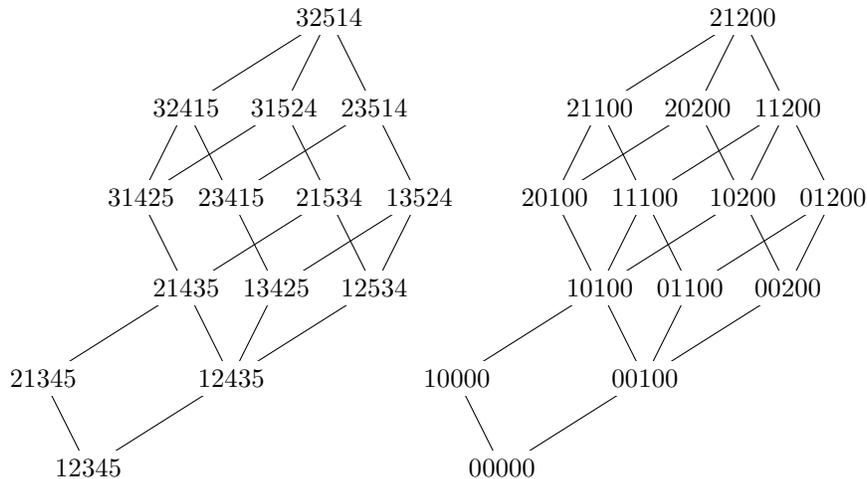
Figure 1: The interval $\Lambda_{32514}$ and its Lehmer codes
\end{center}
The left weak order interval $\Lambda_{32514}$ shown on the left of Figure 1 is not a distributive lattice due to the subinterval $\left[12435, 32415\right]$. Restricted to the Lehmer codes of permutations in $\Lambda_{32514}$, the product order on $\mathbb{N}^5$ refines the left weak order. This is shown on the right of Figure 1. By Theorem~\ref{t:distributivelattice}, this refinement results in a distributive lattice.

Our results relating weak order intervals and distributive lattices are motivated by the existence of nice structure theorems for finite distributive lattices. For example, the fundamental theorem of finite distributive lattices states that any finite distributive lattice is isomorphic to the set $J(P)$ of down-closed subsets of a finite poset $P$, ordered by inclusion. In light of Theorem~\ref{t:distributivelattice}, we construct a finite poset $M_w$ associated to the set of Lehmer codes of permutations in $\Lambda_w$. In Section~\ref{s:baseposet}, we give a chain decomposition of $M_w$ in which the chains are determined by the Lehmer code. The relations between the chains are determined by an extension to the Lehmer code that we introduce in Section~\ref{s:extendedcodes}. The construction of $M_w$ and its properties are summarized by Theorem~\ref{t:codejoinirreducibles}.

Propp \cite{proppdistributive} gave a method for choosing elements uniformly at random from any finite distributive lattice of the form $J(P)$ that uses only the poset $P$. Thus, the description of $M_w$ given in Theorem~\ref{t:codejoinirreducibles} can be combined with this method to choose elements uniformly at random from any weak order interval in $S_n$. 

Our current work is also motivated by questions given at the end of \cite{Weiweakorder} about the rank-generating function of $\Lambda_w$. One question asks which $w \in S_n$ are such that the interval $\Lambda_w$ is rank-symmetric. In Proposition~\ref{p:lowerbound}, we show that there are at least $\left(\lfloor\frac{n}{2}\rfloor\right)!$ such permutations in $S_n$.
\subsection{Preliminaries}
We use the convention that $\mathbb{N} = \{0,1,2,\ldots\}$ and $[n] = \{1,\ldots,n\}$. To specify permutations, we use $1$-line notation. That is, we say $w = w_1w_2 \cdots w_n$ to specify the permutation satisfying $w(i) = w_i$ for all $i \in [n]$.

For any poset $(P,\leq)$, we say that $P$ is ranked if there is a function $\rho:P \rightarrow \mathbb{N}$ satisfying $\rho(x) = 0$ for minimal elements $x \in P$ and $\rho(y) = \rho(x) + 1$ whenever $y$ covers $x$. Whenever $P$ is ranked and finite, the \emph{rank-generating function for $P$} is defined by
\begin{equation*}
F(P,q) = \sum_{x \in P} q^{\rho(x)}.
\end{equation*}
For any poset $(P,\leq)$, a down-closed subset $I \subseteq P$ is called an \emph{order ideal}. That is, a subset $I \subseteq P$ is an order ideal if $y \in I$ whenever $x \in I$ and $y \leq x$. We denote the weak order interval $[\text{id},w]$ by $\Lambda_w$.

For the remainder of this paper, let $n$ be a positive integer.

\vspace{\parskip}
\begin{definition} \label{d:inversiondef}
Let $w \in S_n$ and set
\begin{equation*}
\inv{w} = \{(i,j) \in [n] \times [n] : i < j \text{ and } w(i) > w(j) \}.
\end{equation*}
The set $\inv{w}$ is called the \emph{inversion set of $w$} and each pair $(i,j) \in \inv{w}$ is called an \emph{inversion of $w$}. Regarding $w \in S_n$ as a permutation in $S_{n+1}$ satisfying $w(n+1) = n+1$, set
\begin{equation*}
\ninv{w} = \{ (i,j) \in [n] \times [n+1] : i \leq j \text{ and } w(i) \leq w(j) \}.
\end{equation*}
We call $\ninv{w}$ the \emph{set of non-inversions} of $w$ and each pair $(i,j) \in \ninv{w}$ is called a \emph{non-inversion} of $w$.
\end{definition}
\noindent
The choice to include pairs of the form $(i,i)$ or $(i,n+1)$ in the definition of non-inversion simplifies later characterizations and proofs. Note that $\ninv{w}$ is the complement of $\inv{w}$ relative to the ordered pairs $(i,j) \in [n] \times [n+1]$ satisfying $i \leq j$. In particular, when $(i,j) \in \ninv{w}$, we have $i \leq j$.

\vspace{\parskip}
\begin{definition} \label{d:weakorderlength}
The \emph{length $\ell(w)$} of $w$ is defined by $\ell(w) = \card{\inv{w}}$. The left weak order $(S_n, \bleq)$ is defined as the transitive closure of the relations
\begin{equation*}
v \bleq w \text{ if } w = s_i v \text{ and } \ell(w) = \ell(v) + 1,
\end{equation*}
where $s_i = (i \;\; i + 1)$ is an adjacent transposition in $S_n$. 
\end{definition}
It is known that $(S_n,\bleq)$ is a ranked poset, where length is the rank function.

The right weak order $(S_n,\leq_R)$ has a similar definition where the condition $w = s_i v$ is replaced by $w = v s_i$. Thus $u \leq_R w$ if and only if $u^{-1} \bleq w^{-1}$. The results of our paper can be translated to the right weak order by using the fact that
\begin{equation*}
(\Lambda_w,\leq_R) \cong (\Lambda_{w^{-1}},\bleq).
\end{equation*}
Also, the dual of \cite[Proposition 3.1.6]{bb} states that $[\text{id},wv^{-1}] \cong [v,w]$ for intervals in the left weak order. Thus, our results for principal order ideals can be translated to arbitrary intervals in the left weak order.

For this paper, the following characterization of the left weak order will be more convenient to use than the definition.

\vspace{\parskip}
\begin{lemma} \label{l:weaksubset}
Let $v,w \in S_n$. Then $v \bleq w$ if and only if $\inv{v} \subseteq \inv{w}$. Consequently, we have $v \bleq w$ if and only if $\ninv{w} \subseteq \ninv{v}$.
\end{lemma}
\begin{proof}
This is a dual version of \cite[Proposition 3.1]{permstatslinext}.
\end{proof}
\noindent
For each $i \in [n]$, let $\lehc{i}{w}$ be the number of inversions of $w$ with the first coordinate equal to $i$; that is, $$c_i(w) = \card{{k : (i,k) \in Inv(w)}}.$$ The finite sequence
\begin{equation*}
\leh{w} = (\lehc{1}{w},\ldots,\lehc{n}{w})
\end{equation*}
is called the \emph{Lehmer code for $w$}. 

\vspace{\parskip}
\begin{example}
Let $w = 412563$. The inversions are
\begin{equation*}
(1,2), (1,3), (1,6), (4,6), \text{ and } (5,6).
\end{equation*}
The number of inversions whose first coordinate is $i$ gives the $i$-th coordinate of the Lehmer code. Thus, the Lehmer code of $w$ is $(3,0,0,1,1,0)$.
\end{example}
We view $\bm{c}$ as a function
\begin{equation*}
\bm{c}:S_n \rightarrow \prod_{i=1}^n [0,n-i],
\end{equation*}
mapping each $w \in S_n$ to an $n$-tuple that satisfies the bound $0 \leq \lehc{i}{w} \leq n - i$. It is known (see \cite[Chapter I]{schubertnotes}) that $\bm{c}$ is a bijection and that
\begin{equation*}
\sum_{i=1}^n \lehc{i}{w} = \ell(w).
\end{equation*}
Whenever we need $\lehc{{n+1}}{w}$ to be defined, we make the reasonable convention that $\lehc{{n+1}}{w} = 0$.
\section{Extended codes and the weak order} \label{s:extendedcodes}
We define an extension of the standard Lehmer code. This extended code is used to characterize the weak order in terms of codes and is central to the construction given in Section~\ref{s:baseposet}.

\vspace{\parskip}
\begin{definition} \label{d:extendedcode}
Let $w \in S_n$. For $1 \leq i < j \leq n+1$, define $\lehm{i}{j}{w}$ to be the number of inversions $(i,k) \in \inv{w}$ satisfying $k < j$; that is, $$\lehm{i}{j}{w} = \card{{k<j : (i,k) \in \inv{w}}}.$$ This defines a matrix of values that we call the \emph{extended Lehmer code for $w$}.
\end{definition}
\noindent
The Lehmer code of $w \in S_n$ is easily recovered from the extended Lehmer code of $w$.
\vspace{\parskip}
\begin{lemma} \label{l:extendedextends}
Let $w \in S_n$. Then $\lehc{i}{w} = \lehm{i}{n+1}{w}$ for all $i \in [n]$.
\end{lemma}
\begin{proof}
The number of inversions $(i,k) \in \inv{w}$ satisfying $k < n + 1$ is precisely the number of inversions in $w$ of the form $(i,k)$.
\end{proof}
\begin{example}
Let $w = 31524$. The extended Lehmer code of $w$ (in matrix form) is
\begin{equation*}
\begin{bmatrix}
0 & 0 & 1 & 1 & 2 & 2\\
0 & 0 & 0 & 0 & 0 & 0\\
0 & 0 & 0 & 0 & 1 & 2\\
0 & 0 & 0 & 0 & 0 & 0\\
0 & 0 & 0 & 0 & 0 & 0\\
\end{bmatrix}
\end{equation*}
and the Lehmer code of $w$ is $(2,0,2,0,0)$. The Lehmer code $\leh{w}$ is obtained by reading down the last column of the matrix of $\lehm{i}{j}{w}$.
\end{example}
\noindent
\begin{lemma} \label{l:codeinequality}
Let $v,w \in S_n$ and suppose $v \bleq w$. Then, for all $i \in [n]$ and $j \in [n+1]$, we have
\begin{enumerate}[(a)]
\item $\lehm{i}{j}{v} \leq \lehm{i}{j}{w}$;
\item $\lehc{i}{v} \leq \lehc{i}{w}$.
\end{enumerate}
\end{lemma}
\begin{proof}
Suppose $v \bleq w$. By Lemma~\ref{l:weaksubset}, we have $(i,k) \in \inv{w}$ whenever $(i,k) \in \inv{v}$. Statement (a) follows from Definition~\ref{d:extendedcode}, which, by Lemma~\ref{l:extendedextends}, proves statement (b).
\end{proof}
\begin{remark} \label{r:extendedremark}
There exist $v,w \in S_n$ satisfying the inequality $\lehc{i}{v} \leq \lehc{i}{w}$ for all $i \in [n]$, but $v \not\bleq w$. Thus the code inequality given in Lemma~\ref{l:codeinequality}(b) is not enough to characterize the left weak order. Proposition~\ref{p:weakequivalence} gives an inequality characterization of the left weak order using the extended Lehmer code.
\end{remark}
\noindent
Whether a pair is an inversion or a non-inversion can be detected using the extended Lehmer code. The hypothesis $i \leq j$ below guarantees that either $(i,j) \in \inv{w}$ or $(i,j) \in \ninv{w}$.

\vspace{\parskip}
\begin{lemma} \label{l:inversionequivalence}
Let $w \in S_n$, $i \in [n]$, and $j, k \in [n+1]$. Suppose $i \leq j \leq k$. Then the following are equivalent:
\begin{enumerate}[(a)]
\item $(i,j) \in \ninv{w}$;
\item $\lehm{i}{k}{w} \leq \lehm{i}{j}{w} + \lehm{j}{k}{w}$;
\item $\lehc{i}{w} \leq \lehc{j}{w} + \lehm{i}{j}{w}$.
\end{enumerate}
\end{lemma}
\begin{proof}
By Lemma~\ref{l:extendedextends}, we have $\lehc{i}{w} = \lehm{i}{n+1}{w}$ and $\lehc{j}{w} = \lehm{j}{n+1}{w}$. Thus the specialization $k = n + 1$ proves that (b) $\Rightarrow$ (c). Define the following subsets of $\inv{w}$:
\begin{eqnarray*}
\begin{aligned}
A &= \{(i,l) \in \inv{w} : l < k \};\\
B &= \{(i,l) \in \inv{w} : l < j \};\\
C &= \{(i,l) \in \inv{w} : l = j \};\\
D &= \{(i,l) \in \inv{w} : j < l < k\}.
\end{aligned}
\end{eqnarray*}
It is clear that $A = B \cup C \cup D$ and that the union is pairwise disjoint. By Definition~\ref{d:extendedcode}, we have $\card{A} = \lehm{i}{k}{w}$ and $\card{B} = \lehm{i}{j}{w}$. Therefore
\begin{equation*}
\lehm{i}{k}{w} = \lehm{i}{j}{w} + \card{C} + \card{D}.
\end{equation*}
The remaining implications are proven below by comparing $\card{C} + \card{D}$ to $\lehm{j}{k}{w}$.

Suppose $(i,j) \in \ninv{w}$, so that $w(i) \leq w(j)$ and $|C| = 0$. If $(i,l) \in D$, then $l < k$ and $(j,l) \in \inv{w}$ since $w(j) \geq w(i) > w(l)$. Thus,
\begin{equation*}
\card{C} +\card{D} = \card{D} \leq \lehm{j}{k}{w}.
\end{equation*}
Therefore (a) $\Rightarrow$ (b).

Suppose $(i,j) \in \inv{w}$ so that $w(i) > w(j)$. Suppose that $(j,l) \in \inv{w}$ and that $j < l < k$. Then $(i,l) \in D$ since $w(i) > w(j) > w(l)$. Thus $\card{D} \geq \lehm{j}{k}{w}$. Since $(i,j) \in \inv{w}$, we have $\card{C} = 1$, which implies 
\begin{equation*}
\lehm{i}{k}{w} > \lehm{i}{j}{w} + \lehm{j}{k}{w}. 
\end{equation*}
Specializing to $k = n + 1$ gives the contrapositive of (c) $\Rightarrow$ (a).
\end{proof}
\noindent
The following lemma, which we frequently use in the sequel, is a simple consequence of transitivity on the usual ordering of $\mathbb{N}$.

\vspace{\parskip}
\begin{lemma} \label{l:biconvexity}
Let $w \in S_n$ and let $i, j, k \in [n+1]$. Suppose $i\leq j \leq k$. Then
\begin{enumerate}[(a)]
\item If $(i,j) \in \inv{w}$ and $(j,k) \in \inv{w}$, then $(i,k) \in \inv{w}$;
\item If $(i,j) \in \ninv{w}$ and $(j,k) \in \ninv{w}$, then $(i,k) \in \ninv{w}$.
\item If $(i,j) \in \inv{w}$ and $(i,k) \in \ninv{w}$, then $(j,k) \in \ninv{w}$.
\item If $(i,j) \in \ninv{w}$ and $(i,k) \in \inv{w}$, then $(j,k) \in \inv{w}$.
\end{enumerate}
\end{lemma}
\begin{proof}
Each statement follows from Definition~\ref{d:inversiondef}.
\end{proof}
\noindent
The numerical characterization of the weak order given in Proposition~\ref{p:weakequivalence} below plays a central role in the theorems we obtain. For any pair $(i, j)$, we call the difference $j - i$ the \emph{height} of $(i, j)$.

\vspace{\parskip}
\begin{proposition} \label{p:weakequivalence}
Let $v,w \in S_n$. The following statements are equivalent:
\begin{enumerate}[(a)]
\item The inequality $v \bleq w$ holds in the left weak order;
\item For all $(i,j) \in \ninv{w}$, we have
\begin{equation*}
\lehc{i}{v} \leq \lehc{j}{v} + \lehm{i}{j}{w}.
\end{equation*}
\end{enumerate}
\end{proposition}
\begin{proof}
Suppose $v \bleq w$ and $(i,j) \in \ninv{w}$. By Lemma~\ref{l:codeinequality}, we have
\begin{equation*}
\lehm{i}{j}{v} \leq \lehm{i}{j}{w},
\end{equation*}
and by Lemma~\ref{l:weaksubset}, we have $(i,j) \in \ninv{v}$.  Thus, Lemma~\ref{l:inversionequivalence} implies
\begin{equation*}
\lehc{i}{v} \leq \lehc{j}{v} + \lehm{i}{j}{v}.
\end{equation*}
Combining these inequalities yields $\lehc{i}{v} \leq \lehc{j}{v} + \lehm{i}{j}{w}$. Thus (a) $\Rightarrow$ (b).

For the converse, suppose for a contradiction that $v\not\leq_L w$, thus $\inv{v} \not\subseteq \inv{w}$. Choose a pair $(i,k)$ of minimal height $k - i$, satisfying the property:
\begin{equation*}
\tag{P}
(i,k) \in \inv{v} \tand (i,k) \in \ninv{w}.
\end{equation*}
Lemma~\ref{l:inversionequivalence} implies
\begin{equation*}
\lehc{i}{v} - \lehc{k}{v} > \lehm{i}{k}{v}.
\end{equation*}
By hypothesis, we have $\lehc{i}{v} \leq \lehc{k}{v} + \lehm{i}{k}{w}$ whenever $(i,k) \in \ninv{w}$. Thus,
\begin{equation*}
\lehm{i}{k}{w} \geq \lehc{i}{v} - \lehc{k}{v}.
\end{equation*}
Therefore $\lehm{i}{k}{w} > \lehm{i}{k}{v}$.
Definition~\ref{d:extendedcode} implies the existence of $j < k$ such that $(i,j) \in \inv{w}$ and $(i,j) \in \ninv{v}$.
By Lemma~\ref{l:biconvexity} (c) and (d), we have $(j,k) \in \ninv{w}$ and $(j,k) \in \inv{v}$. Since $k - j < k - i$, this contradicts the minimality of the height of $(i,k)$ with respect to property (P).
\end{proof}
\begin{remark}
Since $\lehc{{n+1}}{v} = 0$ and $\lehm{i}{n+1}{w} = \lehc{i}{w}$, the requirement of Lemma~\ref{l:codeinequality} that $\lehc{i}{v} \leq \lehc{i}{w}$ for all $i \in [n]$ whenever $v \leq_L w$ is contained in Proposition~\ref{p:weakequivalence}.
\end{remark}
\section{The distributive lattice of Lehmer codes for an interval} \label{s:distributive}
We mix partial order and lattice theoretic language in the usual way. When we say ``$(P,\leq)$ is a lattice'' we mean that the join and meet operations are given by the least upper bound and the greatest lower bound, respectively.

By \cite[Section 1.6]{birkhofflatticetheory}, the product space $\mathbb{N}^n$ is a distributive lattice, as is any sublattice of $\mathbb{N}^n$. Thus we use the symbol ``$\leq$'' for the usual order on $\mathbb{N}$, the symbol ``$\pleq$'' for the product order on the product space $\mathbb{N}^n$, and the symbol ``$\bleq$'' for the left weak order on $S_n$. The product order on $\mathbb{N}^n$ is given by
\begin{equation*}
(x_1,\ldots,x_n) \pleq (y_1,\ldots,y_n) \text{ if and only if } x_i \leq y_i \text{ for all } i \in [n].
\end{equation*}
The meet and join on $\mathbb{N}^n$ are given by
\begin{eqnarray*}
\begin{aligned}
(x_1,\ldots,x_n) \vee (y_1,\ldots,y_n) &= (\text{max}\{x_1,y_1\},\ldots,\text{max}\{x_n,y_n\}) \text{ and }\\
(x_1,\ldots,x_n) \wedge (y_1,\ldots,y_n) &= (\text{min}\{x_1,y_1\},\ldots,\text{min}\{x_n,y_n\}).
\end{aligned}
\end{eqnarray*}
\noindent
For an arbitrary $w \in S_n$, consider the subposet $(\leh{\Lambda_w},\pleq)$ of $\mathbb{N}^n$. This is the set of Lehmer codes for all $v \in S_n$ satisfying $v \bleq w$, ordered by the product order $\pleq$. By Lemma~\ref{l:codeinequality}, we know that $v \bleq w$ implies $\bm{c}(v) \pleq \leh{w}$. The converse is false in general, which is shown in the example below. 

\vspace{\parskip}
\begin{example}
Let $w = 32145$ and $w' = 34125$. Then $\leh{w} = (2,1,0,0,0)$ and $\leh{w'} = (2,2,0,0,0)$. It is straightforward to check that $w \not\bleq w'$. By comparing coordinates, we see that $\leh{w} \pleq \leh{w'}$.
\end{example}
The above discussion shows that the set $\leh{\Lambda_w}$ contains as many elements as $\Lambda_w$, but there are more pairs of permutations related by $\pleq$ than by $\bleq$. We use Proposition~\ref{p:weakequivalence} to show that the subset $\leh{\Lambda_w}$ of $\mathbb{N}^n$ is a sublattice of $(\mathbb{N}^n, \pleq)$.

\vspace{\parskip}
\begin{lemma} \label{l:latticeopsclosed}
Let $w \in S_n$. The set $\leh{\Lambda_w}$ of Lehmer codes for the order ideal $\Lambda_w$ is closed under the join and meet of $\mathbb{N}^n$.
\end{lemma}
\begin{proof}
Let $\bm{x},\bm{y} \in \leh{\Lambda_w}$. Let $\bm{x} = (x_1,\ldots,x_n)$ and $\bm{y} = (y_1,\ldots,y_n)$. For some $u_1,u_2 \in S_n$ such that $u_1,u_2 \bleq w$, we have $\bm{x} = \bm{c}(u_1)$ and $\bm{y} = \bm{c}(u_2)$. Let $v \in S_n$ satisfy $\bm{c}(v) = \bm{x} \wedge \bm{y}$. Suppose $(i,j) \in \ninv{w}$.

Suppose, without loss of generality, that $\text{min}\{x_j,y_j\} = x_j$. We have 
\begin{equation*}
x_i \leq x_j + \lehm{i}{j}{w}, 
\end{equation*}
by Proposition~\ref{p:weakequivalence} applied to $u_1$. Since $\text{min}\{x_i,y_i\} \leq x_i$, we have 
\begin{equation*}
\text{min}\{x_i,y_i\} \leq \text{min}\{x_j,y_j\} + \lehm{i}{j}{w}. 
\end{equation*}
Since $\text{min}\{x_i,y_i\} = \lehc{i}{v}$ and $\text{min}\{x_j,y_j\} = \lehc{j}{v}$, it follows that
\begin{equation*}
\lehc{i}{v} \leq \lehc{j}{v} + \lehm{i}{j}{w}.
\end{equation*}
Proposition~\ref{p:weakequivalence} implies $v \bleq w$. Thus $v \in \Lambda_w$. Since $\bm{x} \wedge \bm{y}$ is the Lehmer code for $v$, it follows that $\bm{x} \wedge \bm{y} \in \leh{\Lambda_w}$.

A similar argument proves that $\bm{x} \vee \bm{y} \in \leh{\Lambda_w}$.
\end{proof}
\begin{lemma} \label{l:rankeddistributive}
Every finite distributive lattice is ranked.
\end{lemma}
\begin{proof}
See \cite[Theorem 3.4.1]{ecI} and \cite[Proposition 3.4.4]{ecI}.
\end{proof}
\begin{theorem} \label{t:distributivelattice}
Let $w \in S_n$. The poset $\leh{\Lambda_w}$ is a distributive lattice. Furthermore, we have $F(\Lambda_w,q) = F(\leh{\Lambda_w},q)$.
\end{theorem}
\begin{proof}
Lemma~\ref{l:latticeopsclosed} implies that $\leh{\Lambda_w}$ is a sublattice of $\mathbb{N}^n$. Every sublattice of a distributive lattice is itself distributive, so $\leh{\Lambda_w}$ is a distributive lattice. By Lemma~\ref{l:rankeddistributive}, there is a rank function $\rho$ for $\leh{\Lambda_w}$.

Let $v \bleq w$. Let $\text{id} = v_0 <_L \cdots <_L v_k = v$ be a maximal chain in the weak order interval $[\text{id},v]$. Since $v_{i-1} <_L v_i$, we have $\bm{c}(v_{i-1}) \pleq \bm{c}(v_i)$ by Lemma~\ref{l:codeinequality}. Since $v_i$ covers $v_{i-1}$ in the weak order, we have
\begin{equation*}
\sum_{k=1}^n c_k(v_i) = \ell(v_i) = \ell(v_{i-1}) + 1 = \sum_{k=1}^n c_k(v_{i-1}) + 1.
\end{equation*}
This implies that $\bm{c}(v_i)$ covers $\bm{c}(v_{i-1})$ in the product order. It follows that $\rho(\bm{c}(v_i)) = \rho(\bm{c}(v_{i-1})) + 1$ for $i \in [k]$. Since $\rho(\bm{c}(\text{id})) = \ell(\text{id}) = 0$, we have $\rho(\bm{c}(v)) = \ell(v)$ for all $v \in \leh{\Lambda_w}$. Thus, $\Lambda_w$ and $\bm{c}(\Lambda_w)$ have the same rank-generating function.
\end{proof}
\section{A description of the base poset for $\leh{\Lambda_w}$} \label{s:baseposet}
In this section, fix $w \in S_n$.
\subsection{Identifying the base poset $M_w$}
For any finite poset $P$, we denote the set of order ideals of $P$ by $J(P)$. The set of order ideals of a poset, ordered by inclusion, is a distributive lattice. Conversely, the fundamental theorem of finite distributive lattices states that every finite distributive lattice $L$ is isomorphic to $J(P)$ for some finite poset $P$. We call $P$ the \emph{base poset} for the distributive lattice $L$.

Recall that a \emph{join-irreducible} $z \in L$ is a nonzero lattice element that cannot be written as $x \vee y$, where $x$ and $y$ are nonzero lattice elements. It is known that the base poset $P$ of a distributive lattice $L$ is isomorphic to the set of join-irreducibles for $L$. See \cite[Theorem 3.4.1]{ecI} and \cite[Proposition 3.4.2]{ecI} for details.

In this section, we construct the base poset $M_w$ for $\bm{c}(\Lambda_w)$ by identifying its join-irreducibles.

We denote the $j$-th coordinate of $\bm{x} \in \mathbb{N}^n$ by $\pi_j(\bm{x})$.

\vspace{\parskip}
\begin{definition} \label{d:irreduciblecoordinates}
If $i \in [n]$ and $x \in [\lehc{i}{w}]$, define $\minix{i}{x}{w}$ coordinate-wise by
\begin{equation*}
\pi_j(\minix{i}{x}{w}) = \begin{cases}  0 & \text{if $j < i$;}\\
			  0 & \text{if $(i,j) \in \inv{w}$;}\\
                                  \text{max}\{0,x - \lehm{i}{j}{w}\} & \text{if $(i,j) \in \ninv{w}$.}
                    \end{cases}
\end{equation*}
\end{definition}
\noindent
Note that the coordinates of $\minix{i}{x}{w}$ are as small as possible while satisfying the constraints of Proposition~\ref{p:weakequivalence}. In Proposition~\ref{p:joinirreducible}, we show that the $\minix{i}{x}{w}$ defined in Definition~\ref{d:irreduciblecoordinates} are the join-irreducibles of $\leh{\Lambda_w}$.

\vspace{\parskip}
\begin{example}
Let $w = 3412$. Then $\leh{w} = (2,2,0,0)$ and $\minix{1}{1}{w}$, $\minix{1}{2}{w}$, $\minix{2}{1}{w}$, and $\minix{2}{2}{w}$ are all defined. In general, we have $(i,i) \in \ninv{w}$ and $\lehm{i}{i}{w} = 0$. Thus, Definition~\ref{d:irreduciblecoordinates} implies that the $i$-th coordinate of $\minix{i}{x}{w}$ is always $x$. For $j \neq i$, the $j$-th coordinate is automatically zero unless $j > i$ and $(i,j) \in \ninv{w}$. 

For $i = 2$, there are no such $j$, since $(2,3), (2,4) \in \inv{w}$. Thus, only the second coordinate is nonzero in $\minix{2}{1}{w}$ and $\minix{2}{2}{w}$:
\begin{equation*}
\minix{2}{1}{w} = (0,1,0,0) \text{ and } \minix{2}{2}{w} = (0,2,0,0).
\end{equation*}
For $i = 1$, we have $(1,2) \in \ninv{w}$, but $(1,3),(1,4) \in \inv{w}$. By Definition~\ref{d:irreduciblecoordinates}, we need to find $\text{max}\{0,x - \lehm{1}{2}{w}\}$ to find the second coordinate of $\minix{1}{1}{w}$ and $\minix{1}{2}{w}$. By Definition~\ref{d:extendedcode}, we have $\lehm{1}{2}{w} = 0$.  Thus,
\begin{equation*}
\minix{1}{1}{w} = (1,1,0,0) \text{ and } \minix{1}{2}{w} = (2,2,0,0).
\end{equation*}
\end{example}
\vspace{\parskip}
\begin{lemma} \label{l:ininterval}
Suppose $i \in [n]$ and $x \in [\lehc{i}{w}]$. Then $\minix{i}{x}{w} \in \bm{c}(\Lambda_w)$.
\end{lemma}
\begin{proof}
Let $v \in S_n$ be the permutation such that $\bm{c}(v) = \minix{i}{x}{w}$. We use Proposition~\ref{p:weakequivalence} to show that $v \bleq w$. Thus suppose $(j,k) \in \ninv{w}$. 

There are two cases: either $\lehc{j}{v} = 0$ or $\lehc{j}{v} > 0$.

Suppose $\lehc{j}{v} = 0$. Then $\lehc{j}{v} \leq \lehc{k}{v} + \lehm{j}{k}{w}.$

Suppose instead that $\lehc{j}{v} > 0$. By Definition~\ref{d:irreduciblecoordinates}, we have $(i, j) \in \ninv{w}$ and $\lehc{j}{v} = \lehc{i}{v} - \lehm{i}{j}{w}$. By Lemma~\ref{l:biconvexity}(b), we have $(i, k) \in \ninv{w}$. By Lemma~\ref{l:inversionequivalence}, we have
\begin{equation*}
\lehm{i}{k}{w} - \lehm{i}{j}{w} \leq \lehm{j}{k}{w},
\end{equation*}
and by Definition~\ref{d:irreduciblecoordinates}, we have
\begin{equation*}
\lehc{i}{v} - \lehm{i}{k}{w} \leq \text{max}\{0, \lehc{i}{w} - \lehm{i}{k}{w}\} = \lehc{k}{v}.
\end{equation*}
Adding the inequalities gives
\begin{equation*}
\lehc{i}{v} - \lehm{i}{j}{w} \leq \lehc{k}{v} + \lehm{j}{k}{w}.
\end{equation*}
Since $\lehc{j}{v} = \lehc{i}{v} - \lehm{i}{j}{w}$, it follows that $\lehc{j}{v} \leq \lehc{k}{v} + \lehm{j}{k}{w}$. By Proposition~\ref{p:weakequivalence}, we have $v \bleq w$.
\end{proof}
\begin{lemma} \label{l:uniqueminimal}
Suppose $i \in [n]$ and $x \in [\lehc{i}{w}]$. Then $\minix{i}{x}{w}$ is the unique minimal element of $\bm{c}(\Lambda_w)$ with the $i$-th coordinate equal to $x$.
\end{lemma}
\begin{proof}
We have $(i,i) \in \ninv{w}$ and $\lehm{i}{i}{w} = 0$. Thus,  by Definition~\ref{d:irreduciblecoordinates}, the $i$-th coordinate of $\minix{i}{x}{w}$ is $x$.

Suppose $\bm{y} \in \bm{c}(\Lambda_w)$, satisfying $\pi_i(\bm{y}) = x$. Suppose $(i, j) \in \ninv{w}$. By Proposition~\ref{p:weakequivalence}, we have
$\pi_j(\bm{y}) \geq x - \lehm{i}{j}{w}$. Since $\pi_j(\bm{y}) \geq 0$, we have $\pi_j(\bm{y}) \geq \text{max}\{0, x - \lehm{i}{j}{w}\}$. Therefore, by Definition~\ref{d:irreduciblecoordinates}, each coordinate of $\bm{y}$ is at least as large as the corresponding coordinate of $\minix{i}{x}{w}$.

Uniqueness follows from the finiteness of $\bm{c}(\Lambda_w)$ and the fact that the meet of all elements with the $i$-th coordinate equal to $x$ is an element whose $i$-th coordinate is $x$.
\end{proof}
\begin{lemma}\label{l:minixunique}
Suppose $\minix{i}{x}{w} = \minix{j}{y}{w}$, for some $i,j \in [n]$, $x \in [\lehc{i}{w}]$, and $y \in [\lehc{j}{w}]$. Then $i = j$ and $x = y$.
\end{lemma}
\begin{proof}
Let $v$ be the permutation whose Lehmer code is $\minix{i}{x}{w}$. Since $x > 0$, there is a permutation $u \in \Lambda_w$ such that $u$ is covered by $v$ in the left weak order. The codes of $u$ and $v$ differ in only one coordinate.

Suppose $i \neq j$. Then the $i$-th coordinate or the $j$-th coordinate of $\bm{c}(u)$ is the same as $\bm{c}(v)$. This either contradicts that $\bm{c}(v)$ has the property of being the unique minimal element of $\bm{c}(\Lambda_w)$ with the $i$-th coordinate equal to $x$ or that it is the unique minimal element with the $j$-th coordinate equal to $y$. Thus, $i = j$. Definition~\ref{d:irreduciblecoordinates} then implies that $x = y$.
\end{proof}
\begin{lemma} \label{l:decompose}
Let $\bm{x} = (x_1, \ldots, x_n)$ and suppose $\bm{x} \in \bm{c}(\Lambda_w)$. Then
\begin{equation*}
\bm{x} = \bigvee \minix{i}{x_i}{w},
\end{equation*}
where the join is over all $i \in [n]$ such that $x_i > 0$.
\end{lemma}
\begin{proof}
By Lemma~\ref{l:uniqueminimal}, we have $\minix{i}{x_i}{w} \pleq \bm{x}$ for all $i \in [n]$ such that $x_i > 0$. Therefore,
\begin{equation*}
\bigvee_{i:x_i > 0} \minix{i}{x_i}{w} \pleq \bm{x}.
\end{equation*}
Since the $i$-th coordinate of $\bm{x}$ is $x_i$, the $i$-th coordinate of $\bm{x}$ is $0$ or the same as the $i$-th coordinate of $\minix{i}{x_i}{w}$. Therefore,
\begin{equation*}
\bm{x} \pleq \bigvee_{i:x_i > 0} \minix{i}{x_i}{w}.
\end{equation*}
Combining these inequalities proves the lemma.
\end{proof}
\begin{example}
Let $u = 3214$, $v = 2413$, and $w = 3412$. Then $u \not\in \Lambda_w$ and $v \in \Lambda_w$. We have
\begin{equation*}
\leh{u} = (2,1,0,0) \text{ and } \leh{v} = (1,2,0,0) .
\end{equation*}
By Definition~\ref{d:irreduciblecoordinates}, we have 
\begin{eqnarray*}
\begin{aligned}
\minix{1}{1}{w} &= (1,1,0,0); \\
\minix{1}{2}{w} &= (2,2,0,0); \\
\minix{2}{1}{w} &= (0,1,0,0); \\
\minix{2}{2}{w} &= (0,2,0,0). 
\end{aligned}
\end{eqnarray*}
In each instance, the $i$-th coordinate of $\minix{i}{x}{w}$ is equal to $x$. The additional nonzero coordinates ensure that the requirements of Proposition~\ref{p:weakequivalence} are satisfied. 
 
Note that $\leh{v} = \minix{1}{1}{w} \vee \minix{2}{2}{w}$, but $\leh{u} \neq \minix{1}{2}{w} \vee \minix{2}{1}{w}$. Thus, the hypothesis in Lemma~\ref{l:decompose} that $\bm{x} \in \leh{\Lambda_w}$ is necessary. 
\end{example}
\vspace{\parskip}
\begin{proposition} \label{p:joinirreducible}
The set
\begin{equation*}
M_w = \{ \minix{i}{x}{w} : i \in [n] \tand x \in [\lehc{i}{w}]\}
\end{equation*}
is the set of join-irreducibles for $\bm{c}(\Lambda_w)$.
\end{proposition}
\begin{proof}
Suppose $\bm{y} \vee \bm{z} = \minix{i}{x}{w}$. Then either $\bm{y}$ or $\bm{z}$ has the $i$-th coordinate equal to $x$. Suppose, without loss of generality, that $\bm{y}$ has the $i$-th coordinate equal to $x$. By Lemma~\ref{l:uniqueminimal}, we have $\minix{i}{x}{w} \pleq \bm{y}$. Since $\minix{i}{x}{w}$ is the join of $\bm{y}$ and another element, we also have $\bm{y} \pleq \minix{i}{x}{w}$. Therefore $\minix{i}{x}{w}$ is a join-irreducible of $\bm{c}(\Lambda_w)$.

For the converse, suppose $\bm{y}$ is a join-irreducible of $\bm{c}(\Lambda_w)$. By Lemma~\ref{l:decompose},
\begin{equation*}
\bm{y} = \bigvee_{i : x_i > 0} \minix{i}{x_i}{w}.
\end{equation*}
Since $\bm{y}$ is a join-irreducible, we have $\bm{y} = \minix{i}{x_i}{w}$ for some $i \in [n]$.
\end{proof}
\subsection{A chain decomposition for $M_w$}
We can describe the set $M_w$ defined in Proposition~\ref{p:joinirreducible} more explicitly. There is a partition of $M_w$ into chains.

\vspace{\parskip}
\begin{definition} \label{d:chaindecomposition}
Let
\begin{equation*}
\chn{i}{w} = \{ \minix{i}{x}{w} \in M_w \; : \; 1 \leq x \leq \lehc{i}{w} \},
\end{equation*}
where $\chn{i}{w}$ is possibly empty. We call the sets $\chn{1}{w},\ldots,\chn{n}{w}$ the \emph{chain decomposition} of $M_w$.
\end{definition}
\noindent
The terminology is justified by the following lemma.

\vspace{\parskip}
\begin{lemma}\label{l:chaindecomposition}
Let $\chn{1}{w},\ldots,\chn{n}{w}$ be the chain decomposition of $M_w$. Then each $\chn{i}{w}$ is a chain of $M_w$. Furthermore, we have
\begin{equation*}
M_w = \chn{1}{w} \cup \cdots \cup \chn{n}{w},
\end{equation*}
where the union is pairwise disjoint.
\end{lemma}
\begin{proof}
By Definition~\ref{d:irreduciblecoordinates}, we have $\minix{i}{x}{w} \leq \minix{i}{y}{w}$ whenever $x \leq y$. By Lemma~\ref{l:minixunique}, the chains are pairwise disjoint as sets.
\end{proof}
\begin{lemma} \label{l:lessthan}
Suppose $i < j$ and suppose $\minix{i}{x}{w}, \minix{j}{y}{w}$ are defined. Then
\begin{equation*}
\minix{i}{x}{w} \not\pleq \minix{j}{y}{w}.
\end{equation*}
\end{lemma}
\begin{proof}
By Definition~\ref{d:irreduciblecoordinates}, the $i$-th coordinate of $\minix{i}{x}{w}$ is $x > 0$. Since $i < j$ by hypothesis, the $i$-th coordinate of $\minix{j}{y}{w}$ is $0$. Therefore, we have $\minix{i}{x}{w} \not\pleq \minix{j}{y}{w}$.
\end{proof}
\begin{lemma} \label{l:chains}
Suppose $(i,j) \in \inv{w}$. Then, every element of $\chn{i}{w}$ is incomparable with every element of $\chn{j}{w}$.
\end{lemma}
\begin{proof}
Let $\minix{i}{x}{w} \in \chn{i}{w}$ and let $\minix{j}{y}{w} \in \chn{j}{w}$. If $(i, j) \in \inv{w}$, then by Definition~\ref{d:irreduciblecoordinates}, the $j$-th coordinate of $\minix{i}{x}{w}$ is $0$ and the $j$-th coordinate of $\minix{j}{y}{w}$ is $y > 0$. Therefore, we have $\minix{j}{y}{w} \not\pleq \minix{i}{x}{w}$.

By Lemma~\ref{l:lessthan}, we have $\minix{i}{x}{w} \not\pleq \minix{j}{y}{w}$. Thus, the chains $\chn{i}{w}$ and $\chn{j}{w}$ are pairwise incomparable.
\end{proof}
\begin{lemma} \label{l:ninvrelation}
Suppose $(i,j) \in \ninv{w}$, $x \in [\lehc{i}{w}]$, and $y \in [\lehc{j}{w}]$. Then we have $\minix{j}{y}{w} \pleq \minix{i}{x}{w}$ if and only if $y \leq x - \lehm{i}{j}{w}$.
\end{lemma}
\begin{proof}
If $\minix{j}{y}{w} \pleq \minix{i}{x}{w}$, then by Definition~\ref{d:irreduciblecoordinates}, we have
\begin{equation*}
y \leq \text{max}\{0, x - \lehm{i}{j}{w}\}.
\end{equation*}
Since $y > 0$, we have $y \leq x - \lehm{i}{j}{w}$.

Conversely, suppose that $y \leq x - \lehm{i}{j}{w}$. Then $y \leq \pi_j(\minix{i}{x}{w})$, which implies $\minix{j}{y}{w} \leq \minix{i}{x}{w}$ by Lemma~\ref{l:uniqueminimal}.
\end{proof}
\noindent
The theorem below summarizes important properties of $M_w$. There are no relations between chains $\chn{i}{w}$ and $\chn{j}{w}$ when $(i,j) \in \inv{w}$. Otherwise, if $(i,j) \in \ninv{w}$, then the relations are determined by the extended Lehmer code entry $\lehm{i}{j}{w}$.

\vspace{\parskip}
\begin{theorem} \label{t:codejoinirreducibles}
Let $w \in S_n$ and let
\begin{eqnarray*}
\begin{aligned}
M_w &= \{\minix{i}{x}{w} : i \in [n] \tand x \in [\lehc{i}{w}]\} \text{ and }\\
\chn{i}{w} &= \{\minix{i}{x}{w} : x \in [\lehc{i}{w}]\}.
\end{aligned}
\end{eqnarray*}
\begin{enumerate}[(a)]
\item The set of join-irreducibles for $\leh{\Lambda_w}$ is $M_w$.
\item As distributive lattices, we have $(J(M_w),\subseteq) \cong (\leh{\Lambda_w},\pleq)$.
\item If $i < j$ and $\minix{i}{x}{w}, \minix{j}{y}{w}$ are defined, then $\minix{i}{x}{w} \not\pleq \minix{j}{y}{w}$.
\item If $(i,j) \in \inv{w}$, then every element of $\chn{i}{w}$ is incomparable with every element of $\chn{j}{w}$.
\item If $(i,j) \in \ninv{w}$, $x \in [\lehc{i}{w}]$, and $y \in [\lehc{j}{w}]$, then
\begin{equation*}
\minix{j}{y}{w} \pleq \minix{i}{x}{w} \iff y \leq x - \lehm{i}{j}{w}.
\end{equation*}
\end{enumerate}
\end{theorem}
\begin{proof}
Part (a) is given by Proposition~\ref{p:joinirreducible}. Part (b) can be proved by using \cite[Proposition 3.4.2]{ecI}.

Part (c) is given by Lemma~\ref{l:lessthan}, part (d) is given by Lemma~\ref{l:chains}, and Part (e) is given by Lemma~\ref{l:ninvrelation}.
\end{proof}
\begin{example}
Let $w = 41528637$. Then $\leh{w} = (3,0,2,0,3,1,0,0)$. To construct $M_w$ we first form the chains $\chn{i}{w}$ whenever $\lehc{i}{w} > 0$. Then we add the inter-chain relations using the last part of Theorem~\ref{t:codejoinirreducibles}. To refine the disjoint union of the chains, we need the following values of $\lehm{i}{j}{w}$:  
\begin{equation*}
\lehm{1}{3}{w} = 1, \lehm{1}{5}{w} = 2, \lehm{1}{6}{w} = 2, \lehm{3}{5}{w} = 1, \tand \lehm{3}{6}{w} = 1. \end{equation*}
As $(5,6) \in \inv{w}$, the associated chains are pairwise incomparable.
\begin{center}
\begin{tikzpicture}[scale = 1.1]
        \chain{3}{0}{0.0}{1}
        \chain{2}{1.0}{0.0}{3}
        \chain{3}{2.0}{0.0}{5}
        \filldraw (3.0,1.0) circle (0.05cm)  node[below=5pt] {$\chn{6}{w}$};
        \connectem{2}{5.5}{1.0}{6.5}{0.0}
        \connectem{1}{6.5}{1.0}{8.5}{0.0}
        \connectem{1}{7.5}{0.0}{6.5}{1.0}
        \chain{3}{5.5}{0.0}{1}
        \chain{2}{6.5}{0.0}{3}
        \chain{3}{7.5}{0.0}{5}
        \filldraw (8.5,1.0) circle (0.05cm) node[below=5pt] {$\chn{6}{w}$};
\end{tikzpicture}
\end{center}
\begin{center}
Figure 2: Construction of $M_w$
\end{center} 
We construct the poset $M_w$ in two steps. We begin with the chain decomposition in Definition~\ref{d:chaindecomposition}. Then we use Theorem~\ref{t:codejoinirreducibles}(e) to add relations between the chains. See Figure 2.
\end{example}
\section{Rank-symmetry of $\Lambda_w$}
\noindent
Given a polynomial $f$ with nonzero constant term, we denote by $f^R$ the polynomial
\begin{equation*}
f^R(q) = q^{\text{deg}(f)} f(1/q). 
\end{equation*}
Roughly speaking, this is the polynomial whose coefficients are obtained by reversing the coefficients in $f$. Note that the constant term of $F(\Lambda_w, q)$ is always nonzero.

A polynomial is \emph{symmetric} if the coefficients, when read left to right, are the same as when read right to left. So, a polynomial with nonzero constant term is symmetric if and only if $f = f^R$.
 
A ranked poset $P$ is \emph{rank-symmetric} if its rank-generating function $F(P, q)$ is symmetric. By \cite[Corollary 3.11]{Weiweakorder}, if a permutation $w$ is separable, then the interval $\Lambda_w$ is rank-symmetric. We give another class of rank-symmetric weak order intervals.

Recall that the \emph{dual} $P^*$ of a poset $P$ is a poset on the same set as $P$, such that $x \leq y$ in $P^*$ if and only if $y \leq x$ in $P$. A poset is \emph{self-dual} if $P \cong P^*$. If a ranked poset $P$ is self-dual, then it is rank-symmetric. However, the converse is false. The following proposition is not a characterization of rank-symmetric intervals, but it provides a large class of weak order intervals that are rank-symmetric.
 
\vspace{\parskip}
\begin{proposition} \label{p:selfdual}
Let $w \in S_n$. If $M_w$ is self-dual, then the weak order interval $(\Lambda_w,\bleq)$ is rank-symmetric.
\end{proposition}
\begin{proof}
By Theorem~\ref{t:distributivelattice} and Theorem~\ref{t:codejoinirreducibles}(a), we have 
\begin{equation*}
F(J(M_w),q) = F(\Lambda_w,q). 
\end{equation*}
The result then follows from the fact that $J(P)^* \cong J(P^*)$ for any poset $P$.
\end{proof}
\noindent
There is a standard embedding of $S_m \times S_n$ into $S_{m+n}$: If $v = v_1 \cdots v_m \in S_m$ and $w = w_1 \cdots w_n \in S_n$, then
\begin{equation*}
v \oplus w = v_1 \cdots v_m (w_1 + m)(w_2 + m)\cdots(w_n + m)
\end{equation*}
defines the embedding via $(v,w) \mapsto v \oplus w$. In $S_{m+n}$, each $u \leq_L v \oplus w$ can be decomposed as $v' \oplus w'$, where $v' \bleq v$ and $w' \bleq w$. Therefore, we have
\begin{equation*}
F(\Lambda_{v\oplus w},q) = F(\Lambda_v,q)F(\Lambda_w,q)
\end{equation*}
By \cite[Proposition 3.1.2]{bb}, an alternative characterization of left weak order is given by
\begin{equation*}
u \bleq w \iff \ell(u) + \ell(w u^{-1}) = \ell(w).
\end{equation*}
Using this characterization, it is straightforward to show that 
\begin{equation*}
u \bleq w \iff uw^{-1} \bleq w^{-1} \iff \ell(uw^{-1}) = \ell(w) - \ell(u).
\end{equation*}
It follows that $F(\Lambda_{w^{-1}},q) = F^R(\Lambda_w, q)$.

\vspace{\parskip}
\begin{proposition} \label{p:lowerbound}
For any $w \in S_n$, the interval $\Lambda_{w \oplus w^{-1}}$ is rank-symmetric. It follows that there are at least $\left(\lfloor\frac{n}{2}\rfloor\right)!$ permutations in $S_n$ such that $\Lambda_w$ is rank-symmetric.
\end{proposition}
\begin{proof}
The rank-generating function of $\Lambda_{w\oplus w^{-1}}$ in the left weak order is given by
\begin{eqnarray*}
\begin{aligned}
F(\Lambda_{w\oplus w^{-1}},q) &= F(\Lambda_w,q) F(\Lambda_{w^{-1}},q)\\
&= F(\Lambda_w,q) F^R(\Lambda_w,q).
\end{aligned}
\end{eqnarray*}
Since $(f \cdot f^R)^R = f \cdot f^R$ for any polynomial $f$ with nonzero constant term, it follows that $F(\Lambda_{w \oplus w^{-1}},q)$ is symmetric.
\end{proof}
\section{Counterexamples}
Theorem~\ref{t:distributivelattice} asserts that every weak order interval has a rank-generating function that is the same as the rank-generating function of some distributive lattice. This is not true for arbitrary ranked posets. Thus, it is natural to ask whether the ranked posets similar to weak order intervals in $S_n$ possess this property. For the strong Bruhat order on $S_4$ and the weak order on the Coxeter group $D_4$, we show that there are  intervals that do not have the rank-generating function of a distributive lattice. Thus, Theorem~\ref{t:distributivelattice} does not generalize to the strong Bruhat order or to arbitrary weak order intervals of arbitrary Coxeter groups.

The strong Bruhat order $(S_n, \leq_B)$ is defined similarly to the weak order. The condition $w = s_i v$ where $s_i$ is an adjacent transposition is replaced by the condition $w = t v$ where $t$ is \emph{any} transposition. Under the strong Bruhat order, the lower order ideal of the permutation $w = 3412$ has rank-generating function given by
\begin{equation*}
F((\Lambda_{3412}, \leq_B),q) = 1 + 3q + 5q^2 + 4q^3 + q^4.
\end{equation*}
If there exists a distributive lattice $L$ such that $F(L, q) = F((\Lambda_{3412},\leq_B),q)$, then the dual $L^*$ is a distributive lattice with rank-generating function
\begin{equation*}
F(L^*,q) = 1 + 4q + 5q^2 + 3q^3 + q^4.
\end{equation*}
By the fundamental theorem of finite distributive lattices, there is a finite poset $P$ such that $L^* \cong J(P)$. Such a poset $P$ would have $4$ minimal elements, which means that there would be at least $\binom{4}{2} = 6$ two-element ideals. Thus no such distributive lattice $L$ exists.

The Coxeter group of type $D_4$ has distinguished generating set
\begin{equation*}
S = \{s_1,s_2,s_3,s_4\}
\end{equation*}
subject to the relations
\begin{eqnarray*}
\begin{aligned}
s_i^2 &= 1 \text{ for all } i \in \{1,2,3,4\};\\
(s_is_j)^2 &= 1 \text { for all } i,j \in \{1,3,4\};\\
(s_2s_i)^3 &= 1 \text { for all } i \in \{1,3,4\}.
\end{aligned}
\end{eqnarray*}
Let $w = s_2 s_1 s_3 s_4 s_2 s_4 s_3 s_1 s_2$. This element of $D_4$ appeared in \cite{fbI} as an example of an element with a non-contractible inversion triple. The interval  $(\Lambda_w, \bleq)$ has a rank-generating function given by
\begin{equation*}
F(\Lambda_w,q) = 1 + q + 3q^2 + 3q^3 + 4q^4 + 4q^5 + 3q^6 + 3q^7 + q^8 + q^9.
\end{equation*}
This rank-generating function appears in \cite{stanleylefschetz} in a different context. As stated in that paper, it is straightforward to check that there is no distributive lattice with that rank-generating function.
\begin{center}
ACKNOWLEDGEMENTS
\end{center}
We thank Richard Green for the helpful comments and suggestions. We also thank the referees for their useful and insightful suggestions.
\begin{raggedright}
\bibliographystyle{amsplain}
\bibliography{combinatorics}
\end{raggedright}
\end{document}